\documentclass[12pt,a4]{article}
\usepackage{amsmath}
\usepackage{amssymb}
\usepackage{latexsym}
\usepackage{epsfig}
\usepackage{graphicx}
\usepackage{psfrag}
\usepackage[colorlinks=true, pdfstartview=FitV, linkcolor=blue, citecolor=blue, urlcolor=blue,pagebackref=false]{hyperref}
\hypersetup{urlcolor=blue, citecolor=red}
\usepackage{amsfonts}
\usepackage{amsthm}

\newcommand{\mltext}{}

\newcommand{\radius}{R}

\newtheorem{definition}{Definition}[section]
\newtheorem{theorem}[definition]{Theorem}
\newtheorem{lemma}[definition]{Lemma}

\newcommand{\R}{{\mathbb R}}

\def\tilde{\widetilde}
\def \bfo {\begin {eqnarray*} }
\def \efo {\end {eqnarray*} }
\def \ba {\begin {eqnarray*} }

\def \ea {\end {eqnarray*} }

\def \beq {\begin {eqnarray}}

\def \eeq {\end {eqnarray}}

\def \supp {\hbox{supp}\,}

\def \det {\hbox{det}}

\def \p {\partial}

\def\Z{{\Bbb Z}}

\def\p{\partial}

\def\R{\mathbb R}

\title{Two dimensional invisibility cloaking for Helmholtz equation and non-local boundary conditions }

\author{Matti Lassas\thanks{Institute of Mathematics, Helsinki University of Technology, FIN-02015, Finland.
}\quad and
Ting Zhou\thanks{Department of Mathematics, University of Washington, Seattle, WA 98195, US.
}}

\date{}

\begin{document}

\maketitle

\begin{abstract}

Transformation optics constructions have allowed the design of cloaking devices that steer electromagnetic,  acoustic
and quantum parameters  waves around
a  region without penetrating it, so that  this region is hidden from external observations. The material parameters
used to describe these devices are anisotropic, and singular at the interface between the cloaked and uncloaked regions, making physical realization a challenge. These singular material parameters correspond to singular coefficient functions in
the partial differential equations modeling these constructions and the presence of these singularities
causes various mathematical problems and physical effects on the interface surface.

In this paper, we analyze the two dimensional cloaking for Helmholtz equation when there are sources or sinks present inside the cloaked region. In particular, we consider nonsingular approximate invisibility cloaks based on the truncation of the singular transformations.
Using such truncation we analyze
the limit when the approximate cloaking approaches the ideal cloaking.
We show that,  surprisingly, a non-local boundary condition appears on the inner cloak interface.
This effect in the two dimensional (or cylindrical) invisibility cloaks, which seems to be
 caused by the infinite phase velocity near the interface between the cloaked and uncloaked regions,
is very different to the earlier studied behavior of the solutions in the three dimensional cloaks.


\end{abstract}

\section{Introduction}

There has recently been much activity concerning \emph{cloaking}, or rendering
objects invisible to detection by electromagnetic, acoustic, or other type of waves or physical fields.
Many suggestions to implement cloaking has been based on {\it transformation optics}, that is, designs of
electromagnetic or acoustic devices with customized effects on wave propagation, made possible by taking advantage
of the transformation rules for the material properties of optics.
%
All perfect cloaking devices based on
transformation optics  require anisotropic and singular
 material
parameters, whether the conductivity (electrostatic) \cite{GLU2,GLU3},   index of refraction (Helmholtz)
\cite{Le}, \cite{GKLU1},
permittivity and permeability (Maxwell) \cite{PSS}, \cite{GKLU1}, mass tensor (acoustic)
\cite{GKLU1}, \cite{ChenChan}, \cite{Cummer2}, or
effective mass (Schr\"odinger) \cite{GKLU7,GKLU8,Zhang}.
By {\it singular} material parameters, we mean that at least one of the eigenvalues or the values of the functions
describing the material properties goes to zero or infinity at some points
when the material parameters are represented in Euclidean coordinates,
typically on the interface between the cloaked and uncloaked regions.
Both the anisotropy and singularity  present
serious challenges in trying to physically realize such theoretical plans
using metamaterials. Analogous difficulties are encountered in the study
of invisibility cloaks base on ray-theory \cite{LeoT} and plasmonic resonances \cite{AE,MN}.


%

%
To justify the invisibility cloaking constructions,
one needs to study physically meaningful solutions of
the resulting partial differential equations on the whole domain, including the region where material parameters
become singular. In \cite{GKLU1}, the {\it finite energy solutions} are defined to be at least measurable functions with finite energy in (degenerate) singular weighted Sobolev spaces, and satisfy the equations in distributional sense.

Due to the presence of singular material parameters, or mathematically speaking, partial differential equations
with singular coefficient functions, the question how the waves behave in cloaking devices near the surface where
the material parameters are singular is complicated. Indeed, very different kind of behaviors of solutions are possible:
In the three dimensional case, it is proved in \cite{GKLU1} that the transformation optics construction based
on a blow up map allows cloaking with respect to time-harmonic solutions of the Helmholtz
equation or Maxwell's equations as long
as the object being cloaked is passive. In fact, for the Helmholtz equation,
the object can be an active source or sink. Moreover, in \cite{GKLU1} it is shown that the finite energy solutions
 for the Helmholtz equation
in the three dimensional case  satisfy a hidden boundary condition, namely waves inside the cloaked region
satisfy the  Neumann boundary conditions. For Maxwell's equations, the finite energy solutions inside the cloaked
region need to have vanishing Cauchy data  i.e., the
 hidden boundary conditions are over-determined. This leads to non-existence of finite solutions
 for Maxwell's equations with generic internal currents \cite{GKLU1}.  Physically, this non-existence results is related
 to the so-called extraordinary boundary effects
 on the interface between the cloaked and uncloaked regions \cite{ZCWK}.

Another point of view in dealing with the singular anisotropic design for cloaking devices is to approximate the ideal cloaking parameters by nonsingular, or even nonsingular and isotropic, parameters \cite{GKLU6,GKLU7,GKLU8,KOVW,KSVW,LZ}, which has its advantages in practical fabrication. In the truncation based nonsingular approximate cloaking for three dimensional Helmholtz equation \cite{GKLU8}, when it approaches the ideal cloaking, one can obtain above Neumann hidden boundary condition for the finite energy solution.  In the nonsingular and isotropic approximate cloaking, one can obtain different types of Robin boundary conditions by varying
slightly the way how the approximative cloak in constructed, see \cite{GKLU6,GKLU7,GKLU8}.

Similarly,
for Maxwell's equations it has been studied how the approximate cloak behave on the limit when
the approximate cloaks approach the ideal one \cite{LZ}. We note that for Maxwell's equations there are various
suggestions what kind of limiting cloaks are possible in three dimensions. These suggestions are
 based on constructions where additional layers (e.g.\ perfectly conduction
layer) is attached inside the cloak \cite{GKLU1} or where the ideal cloak corresponds to some of the possible self-adjoint
extensions of Maxwell's equations \cite{Weder,Weder2}.
In  the two dimensional or cylindrical cloaking construction for Maxwell's equations, the eigenvalues of permittivity and the permeability of the cloaking medium
do not only contain eigenvalues approaching to zero  (as in 3D) but also some of the eigenvalues
approach infinity at the cloaking interface.
Then, the electric flux density $D$ and magnetic flux density $B$ may blow up even
when there are no sources inside the cloak and an incident plane wave is scatters from the cloak, see \cite{GKLU3}.
However, if a soft-hard (SH)-surface is included inside the cloak, the solutions behave well.

These above examples show how different behavior the solutions may have, in different type of cloaking devices, near the interface between the cloaked and uncloaked regions.

In this paper, we analyze the two dimensional cloaking for Helmholtz equation when there are sources  present inside the cloaked region.
We start with the nonsingular approximate cloaking based on the truncation of the singular transformation.
Taking the limit when the approximate cloaks approach  the ideal cloak,
we show that a non-local boundary condition appears on the inner cloak interface.
This type of boundary behavior is very different from
that the solutions have in three dimensional case discussed in \cite{GKLU8,ZCWK}.
The main result is formulated as Theorem \ref{thm:main}.
{\mltext We note that in the recent preprint \cite{Nguyen} 
of Hoai-Minh Nguyen 
a different type of formulation, based on a transmission problem,  is given for the non-local boundary condition
appearing in two-dimensional cloaking. In \cite{Nguyen}, also cloaking for more general
second order equations and  quantitative 
convergence properties of the three and two dimensional approximative cloaks 
are analyzed.}

 Physically speaking,
the non-local boundary condition is possible due to the fact that the phase velocity
of the waves in the invisibility cloak approaches infinity near the interface between the cloaked and uncloaked regions,
 even though the group velocity stays finite, see \cite{Chen2}.
We note that as the most important experimental implementations of invisibility cloaks \cite{Sc} have been
based on cylindrical cloaks, the appearance of such boundary condition could also
be studied in the present experimental configurations, at least on micro-wave frequencies.
We also study the eigenvalues, i.e., resonances inside the ideal cloak
corresponding to the non-local boundary condition. As these eigenvalues play an essential role
in the study of almost trapped states \cite{GKLU7} and in the development
of the invisible sensors \cite{AE2,GKLU9}, such resonances can be used to study analogous
constructions in cylindrical geometry.


The rest of the paper is organized as following. In Section \ref{sec:ideal}, we present the basics on transformation optics in the electrostatics setting and apply them to the construction of acoustic ideal cloak for the two dimensional Helmholtz equation. Section \ref{sec:approx} is devoted to the nonsingular approximate acoustic cloaking construction and analysis of behaviors of acoustic waves as it approaches the ideal cloaking.

\section{Perfect acoustic cloaking}\label{sec:ideal}

\subsection{Background: electrostatics}

Our analysis is closely related to the inverse problem for electrostatics,  or Calder\'on's conductivity problem
\cite{AP,C,N1,N2,SyU}.
Let $\Omega\subset \R^d$ be a  domain, at the boundary of which electrostatic measurements are to be made,  and
denote by $\sigma(x)$  the anisotropic  conductivity within.
In the absence of sources, an electrostatic potential $u$ satisfies  a divergence form equation,
\begin{equation}\label{johty}
\nabla\cdot \sigma\nabla u = 0
\end{equation}
on $\Omega$. To uniquely fix the solution $u$ it is enough to give its value, $f$, on the
boundary.  {In the idealized case, one measures,  for all voltage distributions
 $u|_{\p \Omega}=f$ on the boundary the corresponding
current fluxes, $\nu\cdotp \sigma \nabla u$,
where $\nu$ is the exterior unit normal to $\p\Omega$.} Mathematically
this amounts to the knowledge of the Dirichlet--Neumann (DN) map,  $\Lambda_\sigma$.
corresponding to $\sigma$, i.e., the map taking the Dirichlet boundary
values of the solution to \eqref{johty} to the
corresponding Neumann boundary values,
\[
\Lambda_\sigma: \ \ u|_{\p \Omega}\mapsto \nu\cdotp \sigma \nabla u|_{\p \Omega}.
\]
 If $F:\Omega\to \Omega,\quad F=(F^1,\dots,F^d)$, is a
diffeomorphism with $F|_{\p \Omega}=\hbox{Identity}$,
then by making the change of variables $y=F(x)$
and setting $u=v\circ F^{-1}$, we obtain
\begin{equation}\label{johty2}
\nabla\cdot \tilde \sigma\nabla v = 0,
\end{equation}
where $\tilde \sigma=F_*\sigma$ is the push forward of $\sigma$ in $F$,
\beq\label{eqn-transf law}
(F_*\sigma)^{jk}(x)=
\frac 1{\det [\frac {\p F^j}{\p y^k}(y)]}
\sum_{p,q=1}^d \frac {\p F^j}{\p y^p}(y)
\,\frac {\p F^k}{\p y^q}(y)  \sigma^{pq}(y)\Big|_{y=F^{-1}(x)}.
\eeq
This can be used to show
that
\ba
\Lambda_{F_*\sigma}=\Lambda_\sigma.
\ea

Thus, there is a
large
(infinite-dimensional) family of conductivities which all give rise to the
same electrostatic
measurements at the boundary. This observation is due to Luc Tartar (see
\cite{KV2} for an account.)
Calder\'on's inverse problem  for anisotropic
conductivities is then the question of
whether two conductivities with the same DN operator must be
push-forwards of each other. There are a number of positive results in this direction in
two dimensions \cite{ALP,LasU,LeU,Sylv} , but it was shown in \cite{GLU2,GLU3}  in three dimensions
and in \cite{KSVW}  two dimensions that, if one allows  singular maps, then in fact there are counterexamples, i.e., conductivities that are undetectable to electrostatic measurements at the boundary.

From now on, we will restrict ourselves to the two dimensional case.
For each $R>0$, let $B_R=\{|x|\le R\}$ and $\Sigma_R=\{|x|=R\}$ be the central ball and sphere of radius $R$, resp., in $\R^2$, and let ${\it O}=(0,0,0)$ denote  the origin.
To construct an invisibility cloak, for simplicity we use the specific
singular coordinate transformation
$F:\R^2\backslash\{{\it O}\}\to \R^2\backslash B_1,$
given by
 \beq\label{transformation}
x=F(y):=\left\{\begin{array}{cl} y,&\hbox{for } |y|>2,\\
\left(1+\frac {|y|}2\right)\frac{y}{|y|},&\hbox{for }0<|y|\leq 2. \end{array}\right.
\eeq
Letting $\sigma_0=1$ be the homogeneous isotropic conductivity
on $\R^2$, $F$ then defines a conductivity $\sigma$ on
$\R^2\backslash B_1$ by the formula
\beq\label{eqn-transf law2}
\sigma^{jk}(x):=(F_*\sigma_0)^{jk}(x),
\eeq
cf. \eqref{eqn-transf law}.
More explicitly, the matrix $\sigma=[\sigma^{jk}]_{j,k=1}^2$ is of the form
\[
\sigma(x)=\frac{|x|-1}{|x|}\Pi(x)+ \frac{|x|}{|x|-1}(I-\Pi(x)),\quad 1<|x|<2,
\]
where $\Pi(x):\R^3\to \R^3$ is the projection to the radial direction, defined by
\beq\label{proj. formula}
\Pi(x)\,v=\left(v\,\cdotp \frac{x}{|x|}\right)\frac{x}{|x|},
\eeq
i.e., $\Pi(x)$ is represented by the matrix $|x|^{-2}xx^t$, cf.
\cite{KSVW}.

One sees that  $\sigma(x)$ is singular at $\Sigma_1$, that is, in  interface between the cloaked and uncloaked regions,
as one
of its eigenvalues, namely the one corresponding to the radial direction,
tends to $0$, and the other tends to $\infty$ as $|x|$ approaching $1^+$.
We extend  $\sigma$ to $B_1$ as an arbitrary
smooth, nonsingular (bounded from above and below) conductivity there.
Let $\Omega=B_3$; the conductivity $\sigma$ is then a cloaking conductivity on $\Omega$,
as it is indistinguishable from $\sigma_0$, {\it vis-a-vis} electrostatic
boundary measurements of electrostatic potentials (treated rigorously as bounded, distributional solutions of the degenerate elliptic boundary value problem corresponding to $\sigma$ \cite{GLU2,GLU3}).

A similar construction based on a blow up map $F$ was
proposed in Pendry, Schurig and Smith \cite{PSS} to cloak the region $B_1$ in $\R^3$
from observation by electromagnetic waves at a positive frequency;
see also Leonhardt \cite{Le} for a related proposition for the Helmholtz equation in $\R^2$ based
on the use of several leaves of an Riemannian sufrace.

\subsection{Ideal cloaking for the Helmholtz equation with interior sources}

We consider the Helmholtz equation, with source term $p$, of the form
\begin{equation}\label{eqn:Helm}\lambda\nabla\cdot\sigma\nabla u+\omega^2 u=p(x)\qquad \hbox{on }\Omega\end{equation}
corresponding cloaking medium with
the inverse of the anisotropic mass density 
and the bulk modulus given by
\begin{equation}\label{eq:para-perf}\sigma^{jk}=\left\{\begin{array}{cl}\sigma_0^{jk}&\hbox{for }\ |x|>2\\ (F_*\sigma_0)^{jk}&\hbox{for }\ 1<|x|\leq2 \\ \sigma_a^{jk}&\hbox{for }\ |x|\leq 1\end{array}\right.,\quad
\lambda=\left\{\begin{array}{cl}\lambda_0&\hbox{for }\ |x|>2\\ F_*\lambda_0&\hbox{for }\ 1<|x|\leq2 \\ \lambda_a&\hbox{for }\ |x|\leq 1\end{array}\right.
\end{equation}
where $(\sigma_0,\lambda_0)$ corresponds to homogeneous vacuum space and $(\sigma_a,\lambda_a)$ are arbitrary smooth, nondegenerate medium in cloaked region $B_1$.
The push-forward of tensor $F_*\sigma_0$ is defined by \eqref{eqn-transf law} and $F_*\lambda_0$ by
\[
(F_*\lambda_0)(x):=[\mbox{det}(DF)\lambda_0]\circ F^{-1}(x).\]
More specifically, if $(\sigma_0,\lambda_0)=(I, 1)$, one has for $1<|x|<2$,
\[
\sigma(x)=\frac{|x|-1}{|x|}\Pi(x)+ \frac{|x|}{|x|-1}(I-\Pi(x)), \quad
\lambda(x)=\frac{|x|}{4(|x|-1)}\]
are both singular at the cloaking surface $\Sigma:=\Sigma_1=\p
B_1$. 

In the next section, we study the regularized approximate cloaking and obtain the behavior of acoustic waves in a singular medium by taking the limit of waves propagating in the nonsingular medium.

\section{Nonsingular approximate cloaking for the Helmholtz equation with interior sources}\label{sec:approx}

At the present time, we assume
that $\sigma$ and $\lambda$ be homogeneous isotropic in side $B_1$, i.e.,
$(\sigma, \lambda)=(\sigma_a\delta^{jk}, \lambda_a)$ with
$\sigma_a$ and $\lambda_a$ arbitrary positive constants.

To start, let $1<R<2$, $\rho=2(R-1)$ and introduce
 the coordinate transformation
 $F_R:\R^2\backslash B_\rho\to \R^2\backslash B_R$,
 \ba
x:=F_R(y)=\left\{\begin{array}{cl} y,&\hbox{for } |y|>2,\\
\left(1+\frac {|y|}2\right)\frac{y}{|y|},&\hbox{for }\rho<|y|\leq 2.
\end{array}\right.
\ea
We define the corresponding approximate mass tensor $\sigma_R$ and bulk modulus $\lambda_R$ as
\beq \label{eq:para-approx}
\sigma^{jk}_R(x)=\left\{\begin{array}{cl }\sigma^{jk}(x)
&\hbox{for } |x|>R,\\
\sigma_a \delta^{jk} &\hbox{for }|x|\leq R. \end{array}\right. \quad \lambda_R(x)=\left\{\begin{array}{cl }\lambda(x) & \hbox{for } |x|>R,\\
\lambda_a &\hbox{for } |x|\leq R,\end{array}\right.
\eeq
where $\sigma^{jk}$ and $\lambda$ are as in \eqref{eq:para-perf}. Note that then $\sigma^{jk}(x)=
\left(\left(F_{R}\right)_*\sigma_0\right)^{jk}(x)$ and $\lambda(x)=\left(\left(F_R\right)_*\lambda_0\right)(x)$ for $|x|>R$. Observe that, for each $R>1$, the medium is nonsingular,
i.e., is bounded from above and below with, however, the lower bound going to $0$,
and the upper bound going to $\infty$ as $R \rightarrow 1^+$.
Consider the solutions of
\begin{equation}\begin{split}\label{eq:Helm-approx}
(\lambda_R\nabla \cdotp\sigma_R \nabla +\omega^2)u_R=&p\quad\hbox{in }\Omega\\
u_R|_{\p \Omega}=&f,
\end{split}\end{equation}
As $\sigma_R$ and $\lambda_R$ are now non-singular everywhere on $\Omega$,
we have
the standard transmission conditions
on  $\Sigma_R:=\{x:\ |x|=R\}$,
\begin{equation}\label{trans a1}\begin{split}
u_R|_{\Sigma_\radius+}&=u_R|_{\Sigma_\radius-},\\
e_r\cdotp \sigma_R \nabla u_R|_{\Sigma_\radius+}&=
 e_r\cdotp \sigma_R \nabla u_R|_{\Sigma_\radius-},
 \end{split}
\end{equation}
where $e_r$ is the radial unit vector and $\pm$ indicates when the
trace on $\Sigma_R$ is computed as the limit $r\to R^\pm$.

Let $\Omega=B_3$.
Then $u_R$ defines two functions $v_R^\pm$ such that
 \ba
u_R(x)=\left\{\begin{array}{cl} v_R^+(F_R^{-1}(x)),&\hbox{for } R<|x|<3,\\
 v_R^-(x),&\hbox{for } |x|\leq R,\end{array}\right.
\ea
and $v_R^{\pm}$ satisfy
\begin{equation}\label{eq:virtual}\begin{split}
(\nabla^2+\omega^2)v_R^+(y)&=p(F_R(y)) \quad \hbox{in
}\rho<|y|<3,\\ v_R^+|_{\p B_3}&=f,
\end{split}
\end{equation}
and
\beq \label{extra-equation}
(\nabla^2+\kappa^2\omega^2)v_R^-(x)&=&\kappa^2p(x), \quad \hbox{in }|x|<R.
\eeq
where $\kappa^2=(\sigma_a\lambda_a)^{-1}$ is a constant.
Moreover, if we assume $\omega^2$ is not an eigenvalue of the transmission problem, then by the transformation law we have
\[e_r\cdot\sigma_R\nabla u_R\big|_{\partial\Omega}=e_r\cdot\nabla v_R^+\big|_{\p\Omega}.\]
This implies that the DN-map $\Lambda_{\sigma_R,\lambda_R}$ at $\p\Omega$ for the approximate cloaking medium \eqref{eq:para-approx} is the same as the DN-map at $\p\Omega$, denoted by $\Lambda_{\rho}$, of a nearly vacuum domain with a small inclusion present in $B_\rho$.

Next, using polar coordinates $(r,\theta)$, $r=|y|$, and $(\tilde{r},\theta)$, $\tilde r=|x|$, the transmission conditions \eqref{trans a1} on the
surface $\Sigma_R$ yield
\begin{equation} \label{trans a2}\begin{split}
v_R^+(\rho,\theta)&=v_R^-(R,\theta),\\
\rho\, \p_rv_R^+(\rho,\theta)&=\kappa R \, \p_{\tilde r}v_R^-(R,\theta).
\end{split}\end{equation}


We consider the source term $\kappa^2p(x)$ where $p(x)\in C^\infty(\R^2)$ with $\supp p\subset B_{R_0}$ ($0<R_0<1$).
It generates a radiating wave $w(x)\in C^\infty(\R^2)$, namely the solution of
\beq\label{Matti A}
(\nabla^2+\kappa^2\omega^2)w=\kappa^2p\qquad\mbox{in }\;\R^2
\eeq
satisfying the Sommerfeld radiation condition.
Moreover, one can write $w$ as
\begin{equation}\label{eq:w}
w(\tilde r, \theta)=\sum_{n=-\infty}^\infty p_n H_{|n|}^{(1)}(\kappa\omega\tilde r)e^{in\theta},\qquad\mbox{for }\;\tilde r>R_0
\end{equation}
where $H_{|n|}^{(1)}(z)$ and $J_{|n|}(z)$
denote the Hankel and Bessel functions, see \cite{CK}.

In $B_R\backslash\overline{B_{R_0}}$ the function $v_R^-(x)$ differs
from $w$ by an entire solution to the homogeneous equation of \eqref{extra-equation}, and thus for $\tilde r\in(R_0, R)$
\ba
v_R^-(\tilde r,\theta)&=&
\sum_{n=-\infty}^\infty (a_{n}J_{|n|}(\kappa
\omega \tilde r)+p_{n}H^{(1)}_{|n|}(\kappa \omega \tilde r))e^{in\theta},
\ea
with yet undefined $a_{n}=a_{n}(\kappa, \omega; R)$.
Similarly, for  $\rho<r< 3$,
\beq\label{Matti B}
v^+_R(r,\theta)&=&
\sum_{n=-\infty}^\infty\big(c_{n}H^{(1)}_{|n|}(\omega r)
+b_{n}J_{|n|}(\omega r)\big)e^{in\theta},
\eeq
with as yet unspecified $b_{n}=b_{n}(\kappa,\omega; R)$ and
$c_{n}=c_{n}(\kappa, \omega; R)$.

Rewriting the boundary value $f$ on $\p \Omega$ as
\ba
f(\theta)
=\sum_{n=-\infty}^\infty f_{n}e^{in\theta},
\ea
we obtain, together with transmission conditions \eqref{trans a2}, the
following equations for $a_n, \, b_n$ and $c_n$:
\begin{equation}\label{system1}
f_n=b_nJ_{|n|}(3\omega)+c_n H_{|n|}^{(1)}(3\omega),
\end{equation}
\begin{equation}\label{system2}
a_n J_{|n|}(\kappa \omega R)+p_n H^{(1)}_{|n|}(\kappa \omega R)
=b_nJ_{|n|}(\omega
\rho)+c_n H_{|n|}^{(1)}(\omega \rho),\end{equation}
\begin{equation}\label{system3}\begin{split}&\kappa R(\kappa \omega a_n( J_{|n|})'(\kappa \omega R)+\kappa
\omega p_n(H^{(1)}_{|n|})'(\kappa \omega R))\\
&\qquad=\rho(b_n \omega(J_{|n|})'(\omega \rho)+\omega c_n (H_{|n|}^{(1)})'(\omega
\rho)).
\end{split}\end{equation}

Solve for $a_n$ and $c_n$ from
\eqref{system2}-\eqref{system3} in terms of $p_n$ and $b_n$, and
use the  solutions obtained and the equation \eqref{system1}
to solve for $b_n$ in terms of $f_n$ and  $p_n$.
This   yields
\begin{equation}
\label{asymptoticsn}\begin{split}
b_n&=\frac 1{J_{|n|}(3\omega)+s_{n}H_{|n|}^{(1)}(3\omega)}
(f_n+\tilde s_n H_{|n|}^{(1)}(3\omega) p_n),\\
c_n&=s_{n}b_n-\tilde s_{n}p_n,\\
a_n&=t_nb_n-\tilde t_n p_n
\end{split}\end{equation}
where
\[\begin{split}
s_{n}&=\frac{1}{D_n}\left\{\rho J_{|n|}(\kappa\omega R)J_{|n|}'(\omega\rho)-\kappa^2RJ_{|n|}'(\kappa\omega R)J_{|n|}(\omega\rho)\right\},\\
t_n&=\frac{1}{D_n}\left\{\rho H_{|n|}^{(1)}(\omega\rho)J_{|n|}'(\omega\rho)-\rho(H_{|n|}^{(1)})'(\omega\rho)J_{|n|}(\omega\rho)\right\},\\
\tilde s_{n}&=\frac{1}{D_n}\left\{\kappa^2R(H_{|n|}^{(1)})'(\kappa\omega R)J_{|n|}(\kappa\omega R)-\kappa^2 R J_{|n|}'(\kappa\omega R)H_{|n|}^{(1)}(\kappa\omega R)\right\},\\
\tilde t_n&=\frac{1}{D_n}\left\{\kappa^2RH_{|n|}^{(1)}(\omega\rho)(H_{|n|}^{(1)})'(\kappa\omega R)-\rho (H_{|n|}^{(1)})'(\omega\rho)H_{|n|}^{(1)}(\kappa\omega R)\right\}
\end{split}\]
with $D_n$ the common denominator given by
\beq\label{Dn formula}
D_n={\kappa^2RJ_{|n|}'(\kappa\omega R)H_{|n|}^{(1)}(\omega\rho)-\rho J_{|n|}(\kappa\omega R)(H_{|n|}^{(1)})'(\omega\rho)}.\eeq

\subsection{Resonances inside the cloak}

Suppose that the boundary data vanishes, i.e., $f\equiv0$.
Then by \eqref{asymptoticsn}, we have
\[
b_n=\frac{\tilde{s}_n H_{|n|}^{(1)}(3\omega)}{J_{|n|}(3\omega)+s_n H_{|n|}^{(1)}(3\omega)}p_n.
\]
Therefore, one can show
\begin{align}
a_n&=\frac{(t_n\tilde{s}_n-\tilde{t}_ns_n)H_{|n|}^{(1)}(3\omega)-\tilde{t}_nJ_{|n|}(3\omega)}{J_{|n|}(3\omega)+s_n H_{|n|}^{(1)}(3\omega)}p_n \nonumber\\
   &=\frac{\kappa^2 R (H_{|n|}^{(1)})'(\kappa\omega R)l_1-\rho H_{|n|}^{(1)}(\kappa\omega R)l_2}{\rho J_{|n|}(\kappa\omega R)l_2-\kappa^2 R J_{|n|}'(\kappa\omega R)l_1}p_n:=\frac{A_n}{B_n}p_n \label{eq:AnBn}
\end{align}
where
\begin{equation}\label{eq:l1l2}\begin{split}
l_1&=J_{|n|}(\omega\rho)H_{|n|}^{(1)}(3\omega)-H_{|n|}^{(1)}(\omega\rho)J_{|n|}(3\omega),\\
l_2&=J_{|n|}'(\omega\rho)H_{|n|}^{(1)}(3\omega)-(H_{|n|}^{(1)})'(\omega\rho)J_{|n|}(3\omega).
\end{split}\end{equation}


Since for small arguments $0<x\ll1$,
\begin{equation}\label{eq:asympt}\begin{split}J_{|n|}(x)&\sim \frac{1}{n!}\left(\frac{x}{2}\right)^n\quad n\geq1, \quad J_{|n|}'(x)\sim \left\{\begin{array}{ll}-\frac{x}{2}&n=0,\\ \frac{1}{2(n-1)!}\left(\frac{x}{2}\right)^{n-1}&n\geq1,\end{array}\right.\\
H_{|n|}^{(1)}(x)&\sim \left\{\begin{array}{ll}\frac{2i}{\pi}\ln\left(\frac{x}{2}\right)&n=0,\\ -\frac{i(n-1)!}{\pi}\left(\frac{2}{x}\right)^n &n\geq1,\end{array}\right.\quad (H_{|n|}^{(1)})'(x)\sim\left\{\begin{array}{ll}2i x^{-1}/\pi &n=0,\\ {in!2^n}x^{-n-1}/\pi &n\geq1,\end{array}\right.\end{split}\end{equation}
when $n\geq1$, where we denote $f\sim g$ if $f-g=o(g)$ as $x\to 0$, and
\begin{equation}\label{eq:AB}\begin{split}A_n&\sim \frac{i2^n\omega^{-n-1}(n-1)!}{\pi}J_{ n}(3\omega)\big[\omega\kappa^2R(H_{ n}^{(1)})'(\kappa\omega R)+nH_{ n}^{(1)}(\kappa\omega R)\big]\rho^{-n},\\
B_n&\sim \frac{-i2^n\omega^{-n-1}(n-1)!}{\pi}J_{ n}(3\omega)\big[\omega\kappa^2RJ_{ n}'(\kappa\omega R)+nJ_{ n}(\kappa\omega R)\big]\rho^{-n},
\end{split}\end{equation}
for $n\geq 1$.
Now we observe that $|a_n|\to \infty $ as $R\rightarrow 1^+$ if
\begin{equation}\label{eq:cond1}
\big[\omega\kappa^2 R(J_{|n|})'(\kappa\omega R)+|n|J_{|n|}(\kappa\omega R)\big]\bigg|_{R=1}=0.
\end{equation}
Note that then
\begin{equation}\label{eq:cond2}
 \big[\omega\kappa^2R(H_{|n|}^{(1)})'(\kappa\omega R)+|n|H_{|n|}^{(1)}(\kappa\omega R)\big]\bigg|_{R=1}\neq0.
\end{equation}
This implies that if $\kappa$ is outside a discrete set and if $\omega$ is such that \eqref{eq:cond1} and \eqref{eq:cond2} are satisfied by functions $J_{|n|}$ and $H_{|n|}^{(1)}$ for some $n$, then there are sources $p$ for which the $H^1$-norm of the solution $u_R$ goes to infinity in the cloaked region (i.e., when resonance happens) as $R\rightarrow1^+$ (i.e., $\rho\rightarrow0$).
We remark that condition \eqref{eq:cond1} implies condition \eqref{eq:cond2} automatically 
and is equivalent to that the function
\beq\label{Bessel eigenfunction} V_{\pm n}(\tilde r,\theta):=J_{|n|}(\kappa\omega \tilde r)e^{\pm in\theta}\eeq
satisfies the boundary value problem
\begin{equation}\label{eq:BVP}\begin{split}(\Delta+\kappa^2\omega^2)V&=0\qquad \mbox{in }\;B_1,\\
\left[\kappa \tilde r\partial_{\tilde r}V+(-\partial_\theta^2)^{1/2}V\right]\Big|_{\tilde r=1^+}&=0.\end{split}\end{equation}
Next we consider the frequencies $\omega$ for which 
\begin{equation}\label{eq:cond3}
\left\{\begin{split}
\left[\omega\kappa^2 RJ_{|n|}'(\kappa\omega R)+|n|J_{|n|}(\kappa\omega R)\right]\bigg|_{R=1}\neq0,
\\
J_{|n|}(3\omega )\neq0,
\end{split}\right.
\qquad \mbox{for any}\; n\in \Z.
\end{equation}


\subsection{Non-local boundary condition with non-resonant frequencies}

In the following, we show that when we have in
$B_2\backslash\overline{B_R}$, $1<R<2$  the approximative cloaking material parameters, then
for the non-resonant frequencies $\omega$, the boundary condition in \eqref{eq:BVP}
{\mltext holds for all solutions when $R\to 1^+$.}

\begin{lemma}\label{lem:main}
Assume that in $\Omega=B_3$ we have  the material parameters $(\sigma_R,\lambda_R)$.
Moreover, suppose $\omega$ is such that \eqref{eq:cond3} holds.
When $R>1$ is sufficiently close to $1$, then for any source $p\in L^2(\Omega)$ supported compactly in $B_1$ and for $f\in H^{1/2}(\p \Omega)$
 the Helmholtz equation  \eqref{eq:Helm-approx} has a unique solution $u_R$. Moreover, as $R\rightarrow1^+$, Fourier
 coefficients of the solution $v^-_R:=u_R|_{B_R}$ in the cloaked region,
\[
v^-_{R,n}(\tilde r)= \int_{0}^{2\pi} e^{-in\theta}v^-_R(\tilde r,\theta) d\theta
\]
{\mltext the limits $\lim_{R\to 1+}v^-_{R,n}(\tilde r)$ exists and we have
for all $n\in \Z$}
\beq\label{eq:nonlocal}
& &\lim_{R\to 1} ( \kappa \tilde r\partial_{\tilde r}v^-_{R,n}(\tilde r)+|n|v^-_{R,n}(\tilde r))\big|_{\tilde r=R}=0,\\ \label{eq:nonlocal2}
& &\lim_{R\to 1} ( \kappa \tilde r\partial_{\tilde r}v^-_{R,n}(\tilde r)+|n|v^-_{R,n}(\tilde r))\big|_{\tilde r=1}=0.
\eeq

\end{lemma}
\begin{proof}~We suppose that $\omega$ is not an eigenvalue of \eqref{eq:BVP}. Then
\[\begin{split}
v_R^-(\tilde r,\theta)&=\sum_{n=-\infty}^\infty (a_{n}J_{|n|}(\kappa\omega \tilde r)+p_{n}H^{(1)}_{|n|}(\kappa \omega \tilde r))e^{in\theta}\\
&=\sum_{n=-\infty}^\infty\left[\frac{A_n}{B_n}J_{|n|}(\kappa\omega \tilde r)+H_{|n|}^{(1)}(\kappa\omega \tilde r)\right]p_ne^{in\theta}
\end{split}\]
where as in \eqref{eq:AnBn}
\[\begin{split}A_n(R)&=\kappa^2 R (H_{|n|}^{(1)})'(\kappa\omega R)l_1-\rho H_{|n|}^{(1)}(\kappa\omega R)l_2\\
B_n(R)&=\rho J_{|n|}(\kappa\omega R)l_2-\kappa^2 R J_{|n|}'(\kappa\omega R)l_1.
\end{split}\]
{\mltext In following, we sometimes denote shortly $A_n(R)=A_n$ and $B_n(R)=B_n$.}

Denote
\[\Phi_n(\tilde r):=\frac{A_n}{B_n}J_{|n|}(\kappa\omega \tilde r)+H_{|n|}^{(1)}(\kappa\omega \tilde r).\]
Apparently $\tilde{\Phi}(\tilde r, \theta):=\Phi_n(\tilde r)e^{in\theta}$ satisfies the Helmholtz equation
\[(\Delta+\kappa^2\omega^2)\tilde{\Phi}=0\qquad\mbox{in }\,B_R.\]
{\mltext To prove the existence of the limits $\lim_{R\to 1+}v^-_{R,n}(\tilde r)$ and  \eqref{eq:nonlocal}, it is sufficient to show
\begin{equation}\label{eq:bdry}\tilde r\kappa \partial_{\tilde r}\Phi_n(\tilde r)+n\Phi_n(\tilde r)\big|_{\tilde r=R}\rightarrow0\quad\mbox{as }\,R\rightarrow 1^+.\end{equation}
Indeed, by \[\partial_{\tilde r}\Phi_n(\tilde r)=\frac{A_n}{B_n}\kappa\omega J_{|n|}'(\kappa\omega \tilde r)+\kappa\omega (H_{|n|}^{(1)})'(\kappa\omega \tilde r),\]
 and $\lim_{R\to 1}B_n(R)\neq0$ (corresponding to the non-resonant case)}, we have
\begin{equation}\label{eq:bdry_1}\begin{split}&\tilde r\kappa\partial_{\tilde r}\Phi_n(\tilde r)+n\Phi_n(\tilde r)\big|_{\tilde r=R}\\
&=\frac{\kappa^2 R}{B_n}(nl_1+\omega\rho l_2)\left((H_{|n|}^{(1)})'(\kappa\omega R)J_{|n|}(\kappa\omega R)-J_{|n|}'(\kappa\omega R)H_{|n|}^{(1)}(\kappa\omega R)\right)
\end{split}\end{equation}
where $l_1$ and $l_2$ are given by \eqref{eq:l1l2}.

For $n\geq1$, as $\rho\rightarrow0^+$,
\begin{equation}\label{eq:bdry_2}\begin{split}
nl_1+\omega\rho l_2=&\omega\rho \left(J_{|n|-1}(\omega\rho)H_{|n|}^{(1)}(3\omega)-H_{|n|-1}^{(1)}(\omega\rho)J_{|n|}(3\omega)\right)\\
=&\mathcal{O}(\rho^{-n+2}).\end{split}\end{equation}
Combining \eqref{eq:bdry_1}, \eqref{eq:bdry_2} and \eqref{eq:AB}, one has
\begin{equation}\label{eq:cond_order1}\tilde r\kappa\partial_{\tilde r}\Phi_n(\tilde r)+|n|\Phi_n(\tilde r)\big|_{\tilde r=R}=\mathcal{O}(\rho^2) \qquad\mbox{as }\,R\rightarrow1^+~ (\hbox{i.e.}\ \ \rho\rightarrow0^+),\end{equation}
which proves \eqref{eq:bdry}.



For $n=0$, from \eqref{eq:asympt}, one has
\[\begin{split}A_0=&-\frac{2i\kappa^2R}{\pi}(H_0^{(1)})'(\kappa\omega R)J_0(3\omega)\ln\left(\frac{\omega\rho}{2}\right)+\kappa^2R (H_0^{(1)})'(\kappa\omega R)H_0^{(1)}(3\omega)\\
&+\frac{2i}{\pi\omega}H_0^{(1)}(\kappa\omega R)J_0(3\omega)+\mathcal{O}(\rho),\\
B_0 =&\frac{2i\kappa^2 R}{\pi}J_0'(\kappa\omega R)J_0(3\omega)\ln\left(\frac{\omega\rho}{2}\right)-\kappa^2RJ_0'(\kappa\omega R)
H_0^{(1)}(3\omega)\\
&+\frac{2i}{\pi\omega}J_0(\kappa\omega R)J_0(3\omega)+\mathcal{O}(\rho).\end{split}\]
Therefore, 
\beq\label{formula 2}
\partial_{\tilde r}\Phi_0(R)=\frac{A_0}{B_0}\kappa\omega J_0'(\kappa\omega R)+\kappa\omega(H_0^{(1)})'(\kappa\omega R)\eeq
has denominator $B_0$ and numerator
\[
\kappa\omega[A_0J_0'(\kappa\omega R)+B_0(H_0^{(1)})'(\kappa\omega R)] =\frac{2\kappa i}{\pi}W_n(\kappa\omega R)+\mathcal{O}(\rho),\]
where $W_n(x)=H_0^{(1)}(x)J_0'(x)-(H_0^{(1)})'(x)J_0(x)$.
This implies
\[\partial_r\Phi_0(R)\sim \frac{W_n(\kappa\omega R)}
{\kappa R J_0'(\kappa\omega R)\ln\left(\frac{\omega\rho}{2}\right)}\rightarrow0\quad \mbox{as }\rho\rightarrow0^+, \]
i.e., the boundary condition \eqref{eq:bdry} is satisfied for $n=0$ and moreover,
\begin{equation}\label{eq:cond_order2}
r\kappa \partial_r\Phi_0(r)\big|_{r=R}=\mathcal{O}\left(\frac{1}{\ln\left(\frac{\omega\rho}{2}\right)}\right)\qquad\mbox{as }
R\rightarrow1^+.\end{equation}
This proves (\ref{eq:nonlocal}). The equation (\ref{eq:nonlocal2}) follows
similarly by evaluating (\ref{eq:bdry_1}) and (\ref{formula 2})
at $\tilde r=1$ instead of $\tilde r=R$.
\end{proof}

Now we are ready to prove our main result of the limit of the waves $u_R$ of the physical approximate cloaking medium as $R\to 1^+$. We recall that 
$\Sigma=\p B_1$.

\begin{theorem}\label{thm:main}
Let $\omega$ be such that \eqref{eq:cond3} is satisfied.
Assume that $u_R$ is the solution of \eqref{eq:Helm-approx} with $f=0$ and $p\in C^{\infty}_0(B_{R_0})$ with $R_0<1$.
Then as $R\to 1^+$, $u_R$ converges uniformly in compact
subsets of  $B_3\setminus \Sigma$ 
to the limit $u_1$ 
satisfying
\beq\label{eq:u1}
& &(\nabla^2+\kappa^2\omega^2)u_1=\kappa^2p\qquad\mbox{in }\;B_1,\\
\label{eq:u1B}& &\kappa \partial_ru_1+(-\partial_\theta^2)^{1/2}u_1\big|_{\partial B_1}=0,
\eeq
and
\begin{equation}\label{eq:u1-int}
u_1\big|_{B_2\backslash\overline{B_1}}=0.
\end{equation}
\end{theorem}

We note that solutions of \eqref{eq:u1}-\eqref{eq:u1-int} with $f\not=0$ and $p=0$
have been analyzed in \cite{GKLU3}.

\begin{proof} Let $R_0<R_1<R$.
Recall that solution   $w\in C^\infty(\R^2)$ of \eqref{Matti A} is the radiating
solution produced by source $\kappa^2p$ in $\R^2$ and it
 has the expansion \eqref{eq:w} for $\tilde r>R_1$.
{\mltext  Consider the Fourier coefficients
\[
w_n(\tilde r)= \int_{0}^{2\pi} e^{-in\theta}w(\tilde r,\theta) d\theta
\]
and denote $P_n(x):=-(\nabla^2+\kappa^2-n^2)w_n(|x|)$.
As $w\in C^\infty(\R^2)$, we see using integration by parts that
\begin{equation}\label{eqn:est_p}
\|P_n(|x|)\|_{L^2(B_3)}\leq C_M(1+|n|)^{-M} \qquad \mbox{for arbitrary } M>0.
\end{equation}
We consider the Fourier coefficients
\[v^-_{R,n}(\tilde r)=\int_{0}^{2\pi} e^{-in\theta}v^-_R(\tilde r,\theta) d\theta,\quad
v^+_{R,n}(r)=\int_0^{2\pi}e^{-in\theta}v^+_R(r,\theta)d\theta.\]
They satisfy the following problem
\beq\label{new eq 1}
& &(-\nabla^2-\kappa^2\omega^2+n^2)v^-_{R,n}(|x|)=P_n(|x|)\qquad\mbox{for }0\leq|x|\leq R,\\
\label{new eq 2}
& &(-\nabla^2-\omega^2+n^2)v^+_{R,n}(|y|)=0\qquad\mbox{for }\rho\leq|y|\leq2,\\
\label{new eq 3}
& &v^+_{R,n}|_{\partial B_2}=0,\\
\label{new eq 4}
& &v^-_{R,n}|_{\partial B_R^-}=v^+_{R,n}|_{\partial B_\rho^+},\quad \kappa R (v^-_{R,n})'(|x|)|_{\partial B_R^-}=\rho {(v^+_{R,n}})'(|y|)|_{\partial B_\rho^+}.
\eeq
By the transmission conditions \eqref{new eq 4},
we see for  $V^\pm_{R,n}(x)=v^\pm_{R,n}(|x|)$
\[\int_{\partial B_R}\partial_{\tilde r} V^-_{R,n}\overline{V^-_{R,n}}dS_x=\frac{R}{\rho}\int_{\partial B_\rho}\frac{\rho}{\kappa R}\partial_r V^+_{R,n}\overline{V^+_{R,n}}dS_y=\int_{\partial B_\rho}\frac{1}{\kappa}\partial_r V^+_{R,n}\overline{V^+_{R,n}}dS_y.\]
Thus, using integration by parts, we obtain
\[\begin{split}
I_1:=&\frac{1}{\kappa}\int_{B_R}P_n\overline{V^-_{R,n}}dy\\
=&\int_{B_R}(-\nabla^2-\kappa^2\omega^2+n^2)V^-_{R,n}\overline{V^-_{R,n}}dx+\frac{1}{\kappa}\int_{B_2\backslash\overline{B_\rho}}
(-\nabla^2-\omega^2+n^2)V^+_{R,n}\overline{V^+_{R,n}}dy\\
=&-\int_{\partial B_R}\partial_{\tilde r} V^-_{R,n}\overline{V^-_{R,n}}dS_x-\left(\int_{\partial B_2}-\int_{\partial B_\rho}\right)\frac{1}{\kappa}\partial_r V^+_{R,n}\overline{V^+_{R,n}}dS_y\\
&\hspace{-1cm}+\int_{B_R}(|\nabla V^-_{R,n}|^2+(-\kappa^2\omega^2+n^2)|V^-_{R,n}|^2)dx+\frac{1}{\kappa}\int_{B_2\backslash\overline{B_\rho}}(|\nabla V^+_{R,n}|^2+(-\omega^2+n^2)|V^+_{R,n}|^2)dy,
\end{split}\]
and then
\[\begin{split}I_1=&{\int_{B_R}(|\nabla V^-_{R,n}|^2+(-\kappa^2\omega^2+n^2)|V^-_{R,n}|^2)dx+\frac{1}{\kappa}\int_{B_2\backslash\overline{B_\rho}}(|\nabla V^+_{R,n}|^2+(-\omega^2+n^2)|V^+_{R,n}|^2)dy}\\
\geq&\int_{B_R}(-\kappa^2\omega^2+n^2)|V^-_{R,n}|^2dx+\frac{1}{\kappa}\int_{B_2\backslash\overline{B_\rho}}(-\omega^2+n^2)|V^+_{R,n}|^2dy.
\end{split}\]


For $|n|\geq N_0>\max\{\kappa^2\omega^2, \omega^2\}$,
the above and $I_1\leq \|V^-_{R,n}\|_{L^2}\|P_n\|_{L^2}$
 imply first that
\begin{equation}\label{eq:bd}
\left(\|V^-_{R,n}\|_{L^2(B_R)}+\|V^+_{R,n}\|_{L^2(B_2\backslash\overline{B_\rho})}\right)
\leq C_{N_0}\|P_n\|_{L^2(B_3)}
\end{equation}
and second that
\begin{equation}\label{eq:bd2}
\left(\|V^-_{R,n}\|_{H^1(B_R)}+\|V^+_{R,n}\|_{H^1(B_2\backslash\overline{B_\rho})}\right)
\leq C'_{N_0} \|P_n\|_{L^2(B_3)}
\end{equation}
where $C_{N_0}$ and  $C'_{N_0}$ are independent of $R$ and $n$.


By Lemma \ref{lem:main}, for each $n\in\Z$ and $\tilde r\in [0,1)$ and
$r\in (0,2)$ there exists limits
 \beq\label{convergence}
 v_n^-(\tilde r)=\lim_{R\rightarrow1^+}v_{R,n}^-(\tilde r),
\qquad v_n^+(r)=\lim_{R\rightarrow 1^+}v^+_{R,n}(r)\eeq
and we denote  $V^\pm_{n}(x)=v^\pm_{n}(|x|)$.
Let now $0<r_1<R_1<1$. Then by \eqref{eq:bd2} the restrictions $ v_{R,n}^-|_{(r_1,R_1)}$, $R>R_1$
are uniformly bounded in $H^1(r_1,R_1)$. By Sobolev embedding
theorem, the set  $\{v_{R,n}^-|_{[r_1,R_1]};\ R_0<R<1\}$ is relatively compact
in $C([r_1,R_1])\subset H^s(r_1,R_1)$, $1/2<s<1$. Thus any
sequence $(v_{R_j,n}^-|_{[r_1,R_1]})_{j=1}^\infty $ with $R_j\to 1$ has
a  subsequence converging in $C([r_1,R_1])$ which limit  has to coincide with $v_n^-$
by (\ref{convergence}).
Thus $v_{R,n}^-$ have to converge to $v_n^-$ in $C([r_1,R_1])$ and 
hence $V_{R,n}^-$ have to converge to $V_n^-$ in $C(\overline B_{R_1}\setminus  B_{r_1})$ as
$R\to 1$. Similarly, for all $\rho_1>0$ we see using \eqref{eq:bd2}
that  $V_{R,n}^+$ have to converge to $V_n^+$ in $C(\overline B_{2}\setminus B_{\rho_1})$ as $R\to 1$.
Now, as the Sobolev norm $u\mapsto \|u\|_{H^1(B_{R_1}\setminus \overline B_{r_1})}$ is a lower semi-continuous
function in $L^2(B_{R_1}\setminus \overline B_{r_1})$ and \eqref{eq:bd2} holds, we see that
$\|V^-_{n}\|_{H^1(B_{R_1}\setminus \overline B_{r_1})} \leq C'_{N_0} \|P_n\|_{L^2(B_3)}$
for all $0<r_1<R_1<1$. Hence, by (\ref{convergence}) we see using e.g.\ monotone convergence theorem and \cite{KKM} that
\beq\label{final H^1 bound}
\|V^-_{n}\|_{H^1(B_1)}=\|V^-_{n}\|_{H^1(B_1\setminus 0)}=
\lim_{r_1\to 0,R_1\to 1} \|V^-_{n}\|_{H^1(B_{R_1}\setminus \overline B_{r_1})} \leq C'_{N_0} \|P_n\|_{L^2(B_3)}.
\eeq
By \eqref{eqn:est_p} we see that $u_1(\tilde r,\theta)=\sum_{n\in \Z} v^-_n(\tilde r)e^{in\theta}$
is a well defined function in $H^1(B_1)$ satisfying \eqref{eq:u1}.
 By (\ref{eq:bd2}),
\ba
\|V^-_{R,n}\|_{C(\overline B_{R_1}\setminus  B_{r_1})}=\|v^-_{R,n}\|_{C([r_1,R_1])} \leq C_{N_0,r_1,R_1}\|v^-_{R,n}\|_{H^1(r_1,R_1)}
\leq C'_{N_0,r_1,R_1} \|P_n\|_{L^2(B_3)}
,
\ea where $C'_{N_0,r_1,R_1}$ does not depend on $R$ or $n>N_0$.
Thus, using \eqref{eqn:est_p}  
and the convergence of $V_{R,n}^-$  to $V_n^-$ in $C(\overline B_{R_1}\setminus  B_{r_1})$
we see that $u_R$ converge to $u_1$ 
uniformly in compact subsets of $B_1\setminus \{0\}$.
Using equation (\ref{new eq 1}), we see the uniform
convergence in a neighborhood of zero, too.
 Similarly,
we see the uniform convergence in compact subsets of $B_3\setminus 
\overline B_1$.
Moreover, by \cite{LM} and \eqref{eq:u1} the boundary values
$\kappa \partial_ru_1|_{ \p B_1}\in H^{-1/2}(\p B_1)$ and $ u_1|_{\p B_1}\in H^{1/2}(\p B_1)$
are well defined. Now, by (\ref{new eq 1}) and \cite{LM} 
we have $\|\kappa \partial_rV^-_n\|_{H^{-1/2}(\p B_1)}\leq
C_2(1+n^2)\|V^-_{n}\|_{H^1(B_1\setminus \overline B_{R_0})}$ where
 $C_2>0$ is independent of $n$. Thus
using \eqref{eqn:est_p}, \eqref{eq:nonlocal2}, and (\ref{final H^1 bound})
we see that the boundary value $\kappa \partial_ru_1+(-\partial_\theta^2)^{1/2}u_1$,
 vanishes. Hence,
the boundary condition \eqref{eq:u1B} is satisfied.

Equation \eqref{eq:u1-int} for $u_1|_{B_3\backslash \overline{B_1}}$ follows
using  \eqref{eqn:est_p}, \eqref{eq:bd2} and the fact that 
by (\ref{asymptoticsn}), (\ref{Dn formula}) for any fixed $n\in \Z$
we have in \eqref{Matti B}
\ba & & c_n(\rho)=\mathcal{O}(\rho^{|n|}),\quad b_n(\rho)=\mathcal{O}(\rho^{|n|})\qquad\mbox{as }\rho\rightarrow0^+,\quad |n|\geq 1,\\
& &c_0(\rho)=\mathcal{O}((\log \rho)^{-1}),\quad b_0(\rho)=\mathcal{O}((\log \rho)^{-1})\qquad\mbox{as }\rho\rightarrow0^+.\ea

}
\end{proof}

We conclude our discussion by remarking that our analysis also explains the limit of approximate electromagnetic cloaks in the cylindrical case with TE/TM polarized incoming waves, as the solutions of Maxwell's equations in this case satisfy the Helmholtz equation.

\medskip
{\bf Acknowledgements.}
ML was partly supported by Finnish Centre of Excellence in Inverse Problems Research,
Academy of Finland COE 213476 and by Mathematical Sciences Research Institute (MSRI).
TZ was partly supported by MSRI.


\bibliographystyle{amsalpha}

\begin{thebibliography}{A}


%


%
%
%
%
%


\bibitem {AE}
A. Alu and N. Engheta, Achieving transparency with
plasmonic and
metamaterial
coatings, {\it Phys. Rev. E} {\bf 72} (2005),
016623.

\bibitem {AE2}
A. Alu and N. Engheta, Cloaking a Sensor, {\it  Phys. Rev. Lett.} {\bf 102} (2009), 233901.

\bibitem {AP}
K. Astala and L. P\"aiv\"arinta:
Calder\'on's
inverse conductivity problem in the plane. {\it Annals
of Math.},
{\bf 163} (2006), 265-299.

\bibitem {ALP}
K. Astala, M. Lassas and L. P\"aiv\"arinta,
Calder\'on's inverse
problem for anisotropic conductivity in the plane,
{\it Comm. PDE} {\bf 30} (2005), 207--224.

\bibitem{C}
A.P. Calder\'on, On an inverse boundary value
problem, {\it Seminar
on
Numerical
Analysis and
its Applications to
Continuum Physics}, Soc. Brasil.
Mat., R\'io de Janeiro, (1980), 65--73.


\bibitem{ChenChan} H. Chen and C.T. Chan, Acoustic cloaking in three dimensions using acoustic metamaterials, {\it Appl. Phys. Lett.} {\bf 91} (2007), 183518.

	
\bibitem{Chen2}
H. Chen and C. Chan, Time delays and energy transport velocities in three dimensional ideal cloaking,
{\it J. Appl. Phys.} {\bf 104} (2008), 033113.

%
%
\bibitem{CK} D. Colton and R. Kress, {\it Inverse Acoustic and Electromagnetic Scattering Theory},Springer-Verlag, Berlin, 1992.


\bibitem{Cummer2} S. Cummer, et al., Scattering Theory Derivation of a 3D
Acoustic Cloaking Shell,  {\it Phys. Rev. Lett.}  {\bf 100} (2008), 024301.






\bibitem{GKLU1}
A. Greenleaf, Y. Kurylev, M. Lassas and G. Uhlmann, Full-wave invisibility of
active devices at
all frequencies,  {\it Comm.  Math. Phys.}, {\bf 279} (2007), 749.


\bibitem
{GKLU3}
A. Greenleaf, Y. Kurylev, M. Lassas and G. Uhlmann, Improvement of
 cylindrical cloaking with the SHS lining. {\it Opt.
Exp.} {\bf 15} (2007), 12717.



\bibitem{GKLU6}
A. Greenleaf, Y. Kurylev, M. Lassas and G. Uhlmann,
Approximate quantum and acoustic cloaking, To appear in {\it J. Spectral Theory.}
 arXiv:0812.1706v1.

\bibitem{GKLU7} A. Greenleaf, Y. Kurylev, M. Lassas and G. Uhlmann,
Approximate quantum cloaking and almost trapped states,
{\it Phys. Rev. Lett.} {\bf  101} (2008), 220404.


\bibitem{GKLU8}
A. Greenleaf, Y. Kurylev, M. Lassas and G. Uhlmann,
{Isotropic transformation optics:
approximate acoustic and quantum cloaking},
{\it New J. Phys.}, {\bf 10} (2008), 115024.


\bibitem{GKLU9}
A. Greenleaf, Y. Kurylev, M. Lassas and G. Uhlmann, Cloaking vs. shielding in transformation optics, To appear in {\it Phys. Rev. E}, arXiv:0912.1872v1.

\bibitem{GLU2}
A. Greenleaf, M. Lassas and G. Uhlmann, Anisotropic conductivities that
cannot detected by
EIT,
{\it Physiolog. Meas.}  (special issue on
Impedance Tomography),
{\bf 24} (2003), 413.


\bibitem{GLU3}
A. Greenleaf, M. Lassas and G.
Uhlmann, On
nonuniqueness for
Calder\'on's
inverse problem,  {\it
Math. Res. Lett.} {\bf 10} (2003),
685.



%






%

\bibitem{KKM}
T. Kilpel\"ainen, J. Kinnunen, and O. Martio, Sobolev spaces with
zero boundary values on metric spaces. {\it Potential Anal.} {\bf 12} (2000), no. 3,
233Ð247.

\bibitem{KOVW}
R. Kohn, D. Onofrei, M. Vogelius and M. Weinstein, Cloaking via change of variables for the Helmholtz equation,
{\it Comm. Pure Appl. Math.}, {\bf 63} (2010), 0973--1016.


\bibitem{KSVW} R. Kohn, H. Shen, M. Vogelius and M. Weinstein,
Cloaking via change of variables in electrical impedance tomography,
{\it Inverse Prob.} {\bf 24} (2008), 015016.

\bibitem{KV2} R. Kohn and M. Vogelius, Identification of an unknown conductivity by means of measurements at the boundary,
in {\it Inverse problems (New York, 1983)}, SIAM-AMS Proc., {\bf 14}, Amer. Math. Soc., Providence, RI, 1984.





\bibitem{LasU} M. Lassas and G. Uhlmann, Determining Riemannian manifold
from boundary measurements,
{\it Ann. Sci.
\'Ecole Norm. Sup.}, {\bf 34} (2001),
771--787.



\bibitem{Le}
U. Leonhardt,
Optical conformal
mapping, {\it Science} {\bf 312} (2006),
1777.

\bibitem{LeoT}
U. Leonhardt and T. Tyc, Broadband Invisibility by Non-Euclidean Cloaking, {\it Science} {\bf 323} (2009), 110-112.





\bibitem{LeU} J. Lee and G. Uhlmann, Determining
anisotropic real-analytic conductivities by boundary measurements,
{\it Comm. Pure Appl. Math.}, {\bf 42} (1989), 1097--1112.

%
\bibitem{LM}
J.-L.\ Lions and E.\ Magenes, E.
Non-homogeneous boundary value problems and applications. Vol. I., Springer-Verlag, 1972. xvi+357 pp

\bibitem{LZ} H.~Y. Liu and T. Zhou, On approximate electromagnetic cloaking by transformation media, {\it submitted}.



\bibitem{MN} G. Milton and N.-A. Nicorovici, On the
cloaking effects
associated with
anomalous localized resonance, {\it
Proc. Royal Soc. A} {\bf 462}
(2006), 3027--3059.




\bibitem{N1} A. Nachman,
 Reconstructions from boundary measurements.  {\it Ann. of Math.} (2)  {\bf 128} (1988), 531--576.


\bibitem{N2} A. Nachman,
Global uniqueness for a two-dimensional inverse boundary
value problem,
{\it Ann. of Math.} {\bf 143} (1996),  71-96.



\bibitem{Nguyen}
H.-M. Nguyen, Approximate cloaking for the Helmholtz equation via transformation optics and consequences for perfect cloaking, preprint.
$http://cims.nyu.edu/\sim hoaiminh/publication.php$


\bibitem{PSS}
J.B.
Pendry, D. Schurig, and D.R. Smith, Controlling Electromagnetic
Fields,
{\it Science}  {\bf 312} (2006), 1780.



\bibitem{Sc} D. Schurig, J. Mock, B. Justice, S. Cummer,
J. Pendry, A. Starr and
D. Smith, Metamaterial electromagnetic cloak at microwave frequencies,
{\it Science} {\bf 314} (2006), 977.



\bibitem{Sylv} J. Sylvester,
An anisotropic inverse
boundary value problem, {\it Comm. Pure Appl.
Math.} {\bf 43} (1990),
201--232.



\bibitem{SyU} J. Sylvester and G. Uhlmann, A global uniqueness theorem for
an
inverse boundary value problem, {\it Ann. Math.} {\bf 125} (1987), 153.


\bibitem{Weder} R. Weder, A rigorous analysis of
high-order electromagnetic invisibility cloaks,
{\it J. Phys. A: Math. Theor.}, {\bf  41} (2008), 065207.

\bibitem{Weder2} R. Weder, The boundary conditions for electromagnetic
invisibility cloaks, {\it J. Phys. A: Math. Theor.} {\bf 41} (2008), 415401.




\bibitem{Zhang} S. Zhang, D. Genov, C. Sun and X. Zhang, Cloaking of matter
waves,
{\it Phys. Rev. Lett.} {\bf 100} (2008), 123002.


\bibitem{ZCWK}
B. Zhang, H. Chen, B.-I. Wu and J. A. Kong,
{Extraordinary surface voltage effect in the invisibility cloak with an active device inside},
{\it Phys. Rev. Lett.}, {\bf 100} (2008), 063904.

%


\end{thebibliography}

\end{document}